\documentclass[11pt, reqno]{amsart}

\usepackage{amssymb,amsthm}
\usepackage[T1]{fontenc}
\usepackage{pdfpages}
\usepackage{enumitem,footnote,color}

\usepackage{caption}
\usepackage{tikz}
\usepackage{dsfont}

\usepackage[colorlinks=True]{hyperref}

\evensidemargin\oddsidemargin

\newtheorem{thm}{Theorem}%[section]
\newtheorem{step}{Step}%[section]
\newtheorem{lem}[thm]{Lemma}

\newtheorem{cor}[thm]{Corollary}

\DeclareMathOperator{\vol}{Vol}

\newcommand{\si}{\sigma}
\newcommand{\om}{\omega}
\newcommand{\Si}{\Sigma}
\newcommand{\X}{\mathbb{S}}
\newcommand{\N}{\mathbb{N}}
\newcommand{\R}{\mathbb{R}}

\newcommand{\la}{\lambda}
\newcommand{\al}{\alpha}

\renewcommand{\leq}{\leqslant}
\renewcommand{\geq}{\geqslant}
\renewcommand{\epsilon}{\varepsilon}

\begin{document}

\baselineskip=16pt

\title{On the first hitting time of a high-dimensional orthant}

\date{\today}

\author{Emmanuel Humbert}

\address{Université de Tours, Institut Denis Poisson, Tours, France. {\em Email}: {\tt emmanuel.humbert@univ-tours.fr}}

\author{Kilian Raschel} 

\address{CNRS, International Research Laboratory France-Vietnam in mathematics and its applications, Vietnam Institute for Advanced Study in Mathematics, Hanoï, Vietnam. {\em Email}: {\tt raschel@math.cnrs.fr}}

\thanks{EH is supported by the project Einstein-PPF (\href{https://anr.fr/Project-ANR-23-CE40-0010}{ANR-23-CE40-0010}), funded by the French National Research Agency. KR is supported by the project RAWABRANCH (\href{https://anr.fr/Project-ANR-23-CE40-0008}{ANR-23-CE40-0008}), funded by the French National Research Agency.}



\begin{abstract} 
We consider a collection of independent standard Brownian particles (or random walks), starting from a configuration where at least one particle is positive, and study the first time they all become negative. This is clearly equivalent to studying the first hitting time from the negative orthant or the first exit time from the complement of the negative orthant.  While it turns out to be possible to compute the distribution of these hitting times for one and two particles, the distribution (and even its tail asymptotics)\ is not known in closed form for three or more particles. In this paper we study the tail asymptotics of the distribution as the number of particles tends to infinity. Our main techniques come from spectral geometry: we prove new asymptotic estimates for the principal eigenvalue of the complement of a high-dimensional orthant, which we believe are of independent interest.
\end{abstract}

\maketitle

\section{Introduction and main results}
\label{sec:introduction}

Consider $d\geq 1$ standard Brownian particles $B_1(t),\ldots,B_d(t)$ (with zero drift and unit variance)\ which at $t=0$ are not all non-positive (i.e., at least one of them starts with a positive value). Our main quantity of interest is the first hitting time of the negative orthant, or equivalently the first exit time from the complement of the negative orthant:
\begin{equation}
\label{eq:def_tau_d}
    \tau_d =\inf\{ t>0 : B_1(t)<0,\ldots, B_d(t)<0\}.
\end{equation}
In the probabilistic literature, this problem can be interpreted as a capture problem \cite{Bramson1991} or an absorption problem \cite{Ke-91}. This first passage time problem has also attracted some attention in the physics community, see \cite{BNKr-10a,BNKr-10b,BrMaSc-13}, and is inspired by classical experiments in thermodynamics. In fact, if we interpret the positive half-line as an almost confined space with a small hole through which particles can escape, we can imagine the following experiment. Let $d$ particles start from the confined space, with a ballistic motion to the left. Then the particles eventually all leave the space and go to the negative half-line. If the same particles follow Brownian motions, how likely is it that they will all have left the room at a given time, especially in the regime of infinitely many particles? 

To formulate the problem precisely, let $x\in\mathbb R^d$ be an arbitrary (starting)\ point and let $\mathbb P_x(\cdot)$ be the probability conditional on the event $x=(B_1(0),\ldots,B_d(0))$. It is proved in \cite[Cor.~1.3]{DB-86} and in \cite[Cor.~1]{BaSm-97} that the survival probability $\mathbb P_x(\tau_d>t)$ behaves asymptotically when $t\to\infty$ as
\begin{equation}
\label{eq:survival_probability}
    \mathbb P_x(\tau_d>t) \sim \frac{V_d(x)}{t^{p_d/2}},
\end{equation}
where $V_d(x)$ is the unique (up to multiplicative scaling)\ harmonic function in the cone $C=\mathbb R^d\setminus \mathbb R_-^d$ (the complement of the negative orthant)\ and the exponent $p_d$ is given by
\begin{equation}
\label{eq:def_pd}
    p_d=\sqrt{\lambda_1 +\bigl(\tfrac{d}{2}-1\bigr)^2}-\bigl(\tfrac{d}{2}-1\bigr),
\end{equation}
with $\lambda_1=\lambda_1(d)$ defined as follows. It is (see e.g.\ \cite[Eq.~(1.3)]{BaSm-97})\ the smallest (positive)\ eigenvalue $\lambda$ of the Dirichlet problem for the Laplace-Beltrami operator $\Delta_{\mathbb S^{d-1}}$ on the sphere $\mathbb S^{d-1}\subset\mathbb R^d$
\begin{equation}
\label{eq:Dirichlet_problem}
     \left\{
\begin{array}{rll}
     -\Delta_{\mathbb S^{d-1}}m&=\ \lambda m & \text{in } C,\\
     m&=\ 0& \text{in } \partial C,
     \end{array}
     \right.
\end{equation}
where $C$ is, as above, the complement of the negative orthant. Remember that in general the quantity $\lambda_1$ cannot be computed in closed form. An asymptotic statement similar to \eqref{eq:survival_probability} is proved in the discrete setting (random walks)\ in \cite[Thm~1]{DeWa-15}.

\begin{figure}
    \centering
    \includegraphics[height=3.5cm]{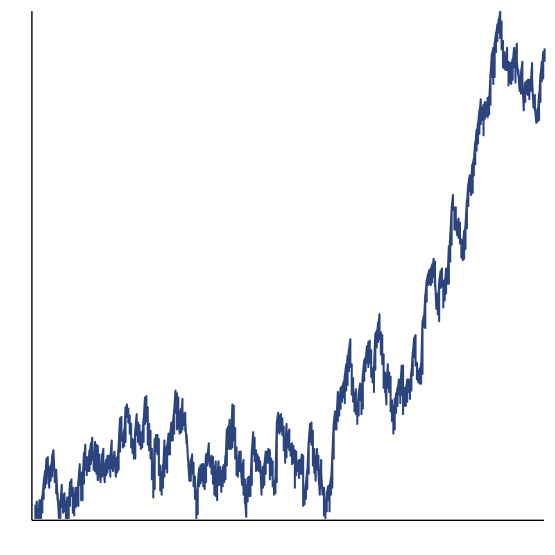}
    \quad\quad
    \includegraphics[height=3.5cm]{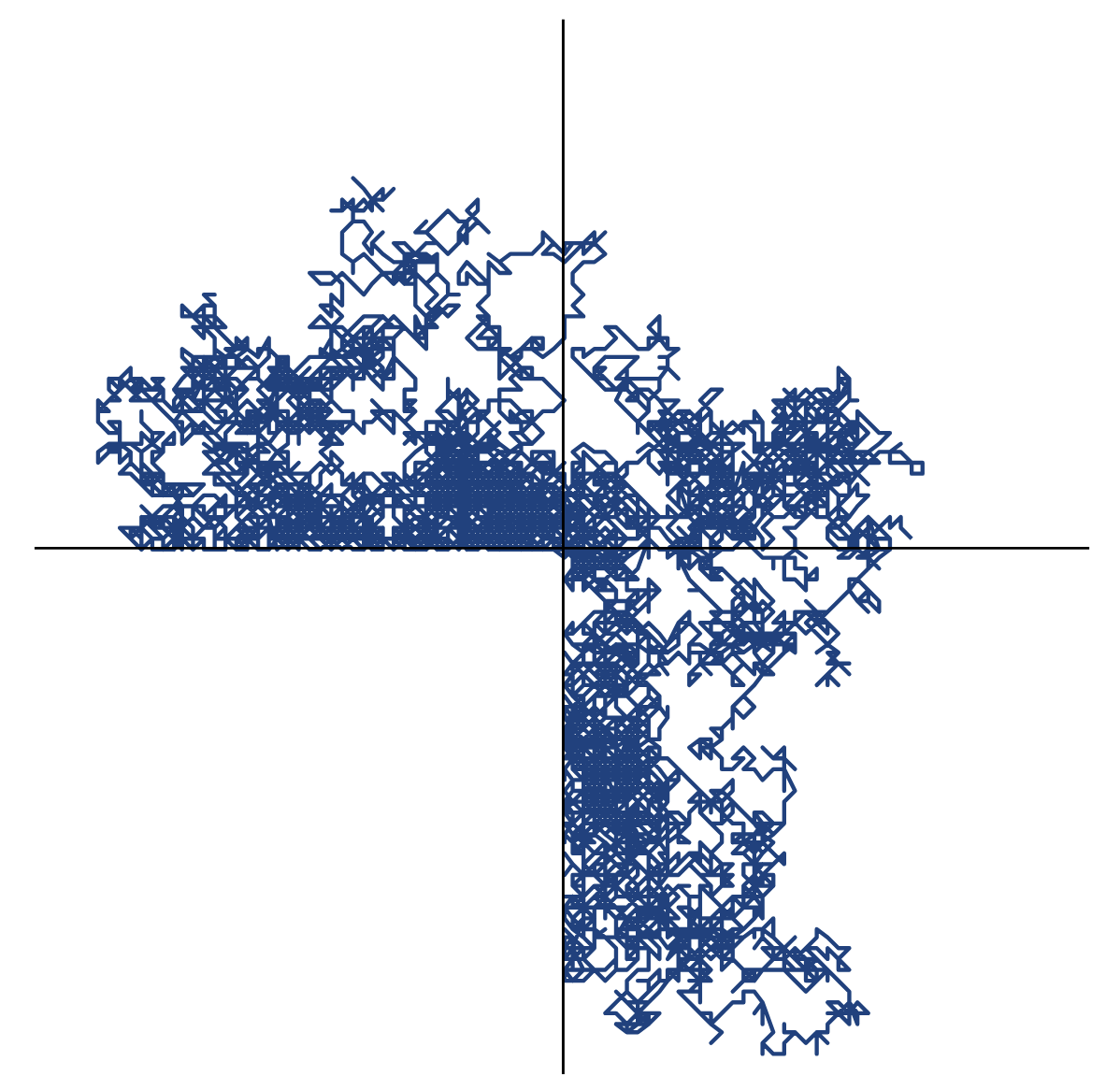}
    \quad\quad
    \includegraphics[height=3.5cm]{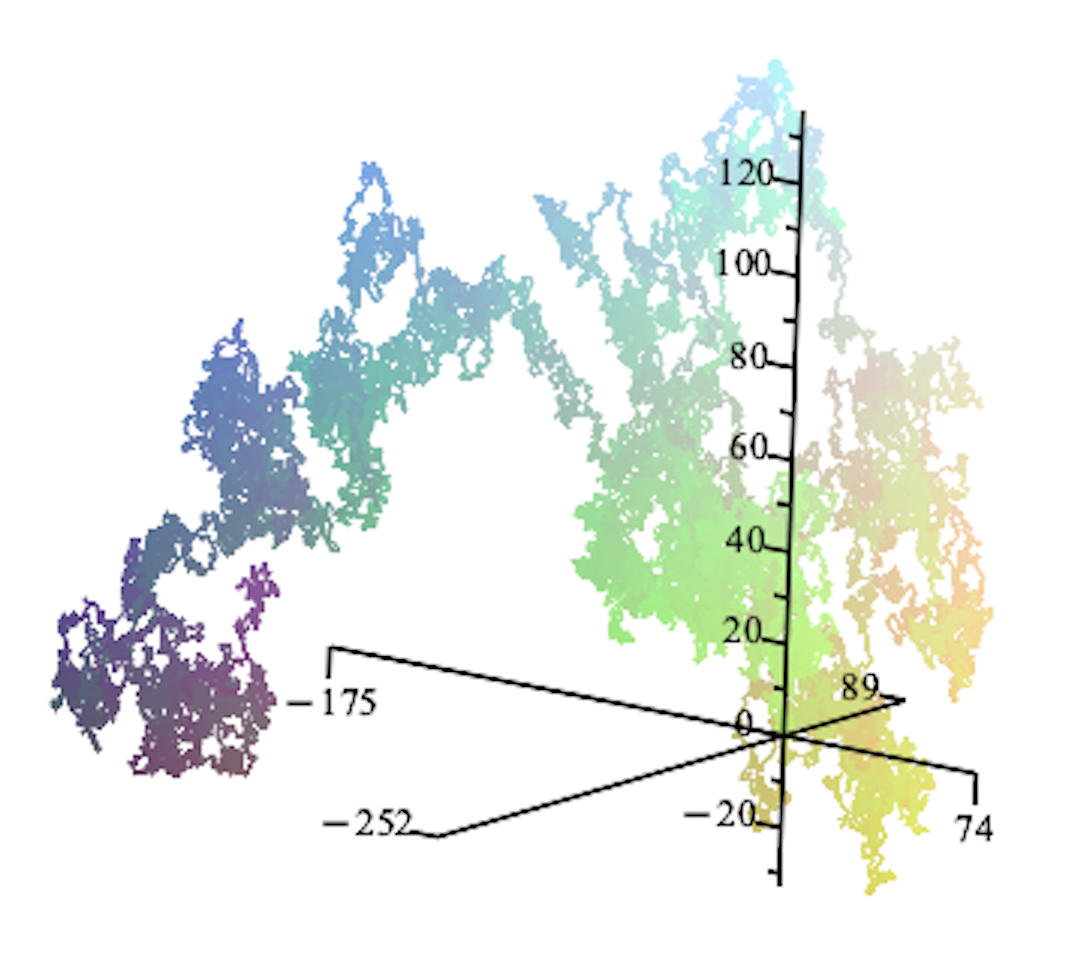}
    \caption{Random walks staying in the complement of the $d$-dimensional orthant, for $d=1,2,3$}
    \label{fig:123}
\end{figure}

It is well known that $p_1=1$ (in dimension $1$, $\tau_1$ in \eqref{eq:def_tau_d} is the first passage time at $0$ of a standard Brownian motion started at a positive value, whose distribution is classical). See Figure~\ref{fig:123}. It is also known that $p_2=\frac{2}{3}$. In dimension $2$, we should mention that the enumeration of walks in the complement of a quadrant has recently attracted some attention, see \cite{BM-16,Mu-19,BM-23,BMWa-24}. In particular, the exponent $\frac{p_2}{2}=\frac{1}{3}$ is obtained in \cite{BM-16}, see the remarks following theorems~1 and 4 in that paper; see also \cite[Thm~3.1]{Mu-19} and \cite[Cor.~4.3]{BMWa-24}. As pointed out in  \cite[Sec.~4]{Mu-19}, the value of $p_3$ (and a fortiori of $p_d$ for $d>3$)\ is unknown. However, the approximation $p_3=0.4542$ is derived in \cite[Sec.~7.1]{BoPeRaTr-20}; see also \cite[Tab.~1]{BNKr-10b}. It would be interesting to certify the previous approximation, as was done in \cite{DaSa-20} for linear transformations of three-dimensional orthants.

Arguing on the number of particles and on \eqref{eq:survival_probability}, it is intuitively clear that $p_d$ is decreasing and goes to $0$ as $d\to\infty$. In this paper we explore the asymptotic behavior of $p_d$, or equivalently of $\lambda_1$ in the high dimension regime.

\begin{thm}
\label{thm:main_estimate}
We have $\lim_{d\to\infty} \frac{\log \lambda_1}{d} = \log ( \frac{1}{2} )$. More precisely, there exist constants $c,C>0$ such that for all $d\geq 1$, $\lambda_1\in \bigl[\frac{cd}{2^d},\frac{Cd^3}{2^d}\bigr]$.
\end{thm}
Obviously, the quantity $\frac{1}{2^d}$ can be interpreted as the fraction of the volume of the positive orthant. However, to our knowledge there is no easy way to translate this elementary observation into a direct estimate of $\lambda_1$. As an immediate consequence of Theorem~\ref{thm:main_estimate}, we have:
\begin{cor}
The critical exponent $p_d$ of the survival probability \eqref{eq:survival_probability} goes to $0$ and is equivalent to
\begin{equation*}
    p_d\sim \frac{\lambda_1}{d}.
\end{equation*}
\end{cor}

Let us conclude this introduction by mentioning five challenging open problems. First, it would be interesting to study the joint asymptotics of the survival probability $\mathbb P_x(\tau_d>t)$ when $d,t\to\infty$  simultaneously. This seems to be a difficult problem, since the asymptotics \eqref{eq:survival_probability} only holds in a regime where $d$ is fixed. However, there are explicit expressions for the survival probability as infinite sums of eigenfunctions, see \cite[Thm~1]{BaSm-97}. Note that when $t=\frac{T}{d}$, the joint asymptotics of the survival probability when $d,T\to\infty$ is studied in \cite{BNKr-10b} and leads to interesting conjectures based on an approximation of cones with spherical caps \cite{BNKr-10a}.

As a second problem, we could consider Brownian motions $B_1(t),\ldots,B_d(t)$ with their own variance $\sigma_1^2,\ldots,\sigma_d^2$ and ask the same asymptotic question as in \eqref{eq:survival_probability}. Technically, this reduces to estimating the principal eigenvalue $\lambda_1$ of the complement of the linear transformation of an orthant. The quantity $\lambda_1$ should also decay exponentially fast to zero.

We could even tackle a more general problem by adding a drift to each Brownian particle. In this case, the problem becomes even more complicated because the associated multidimensional random walk has a drift and, to the best of our knowledge, there is no general expression for the asymptotics of the survival probability in the case of a random walk with drift, except for partial results \cite{Du-14,GaRa-16}.

Moreover, there is a deep connection between the exponent $p_d$ in \eqref{eq:def_pd} and the following occupation time problem. Considering a Brownian motion $B(t) = (B_1(t),\ldots,B_d(t))$ starting from the origin, one can define the occupation time of the orthant $\mathbb R_+^d$ as
\begin{equation*}
    T=\int_0^1 \mathds{1}_{\{B(t)\in \mathbb R_+^d\}}dt.
\end{equation*}
It is well known that in dimension $d=1$, $T$ follows the arcsine law. On the other hand, the distribution of $T$ is unknown as soon as $d\geq 2$. However, it was proved in \cite[Thm~2]{MeWe-95} that as $t\to0$, $\mathbb P(T\leq t)/t^{p_d/2}$ is bounded from below and above by two positive constants. (This result actually holds for a much larger family of cones.)\ In this context,  Theorem~\ref{thm:main_estimate} provides a precise approximation to the short-time behaviour of the occupation time of a large-dimensional orthant.

Finally, in a broader setting, this work initiates the study of random walks in cones in the high-dimensional regime. When there is no cone constraint, high-dimensional random walks are quite well understood, see e.g.\ \cite{KaMa-22} for a study of a random walk when time and dimension simultaneously go to infinity. It would be interesting to understand how a high-dimensional cone structure affects this joint limit behaviour. 

\section{Proof of Theorem~\ref{thm:main_estimate}}

In the whole proof, if $d \in \N^*=\{1,2,3,\ldots\}$ and $A \subset \X^d$, $| A |_d$ will denote the volume of $A$ with respect to the standard metric $\si_d$ of $\X^d$. In particular, $| \mathbb S^d |_d$ is denoted by $\omega_d$. We will also set  \begin{equation}
    \label{eq:def_domain_Ud}
    U_d = \bigl(\mathbb R^d\setminus \mathbb R_-^d\bigr)\cap \mathbb S^{d-1}.    
\end{equation}
Theorem~\ref{thm:main_estimate} aims to estimate the first Dirichlet eigenvalue $\la_1(d)$ of $U_d$ for the Dirichlet problem~\eqref{eq:Dirichlet_problem}.

\subsection{A preliminary lemma} 
For the upper bound part of the proof, we will need the following crucial lemma, to be proved in  Section~\ref{prooflem}.
\begin{lem}
\label{volsi}
Define for $d \in \N^*$ and $a \in [0,1]$
\begin{equation}
\label{eq:def_dom}
    \Si_d(a) = [-a,1]^d \cap \X^{d-1}.
\end{equation}
Let $(a_d)_{d\geq 1}$ be any sequence such that $a_d \ll d^{-\frac{3}{2}}$. 
Then there exists $C>0$ such that for all $d\geq 1$,
\begin{equation*}
   \bigl\vert\Si_d(a_d) \bigr\vert_{d-1} \leq C \frac{\om_{d-1}}{2^d}.
   \end{equation*}
\end{lem}
To the best of our knowledge, the volume of $\Si_d(a)$ does not admit a closed-form expression.

\subsection{The lower bound}
  We denote by $H^1(\X^d)$ the Sobolev space of $L^2$-functions on $\X^d$ whose all  derivatives of order $1$ are also in $L^2$. For any $u \in H^1(\X^d)$, $u \not= 0$, define 
  \begin{equation*}
      Q(u) = \frac{ \frac{4(d-1)}{d-2} \int_{\X^d} |\nabla u|^2 d\si_d + d(d-1) \int_{\X^d} u^2 d\si_d}{  \left( \int_{\X^d} |u|^{\frac{2d}{d-2}} d\si_d \right)^{\frac{d-2}{d}} }.
  \end{equation*}.
  The function $Q$ is called the Yamabe functional and appeared naturally in the solution of the famous Yamabe problem (see e.g.\ \cite{LP-87} for more information). 

It is well known (see for instance \cite[Thm~3.2]{LP-87}) that the infimum of $Q$ over the set $H^1(\X^d)\setminus \{0\}$ is attained for the constant function $u=1$. This  implies that for all $u \in H^1(\X^d)$, $u\not=0$, $Q(u) \geq Q(1)$ and thus 
\begin{equation} 
\label{sob_ineq}
  \left( \int_{\X^d} |u|^{\frac{2d}{d-2}} d\si_d \right)^{\frac{d-2}{d}} \leq \frac{1}{d(d-1) \om_d^{\frac{2}{d}}} \left( \frac{4(d-1)}{d-2} \int_{\X^d} |\nabla u|^2 d\si_d + d(d-1) \int_{\X^d} u^2 d\si_d\right).
\end{equation}
 Let $u_{d+1}$ be an eigenfunction of 
$U_{d+1}$, see \eqref{eq:def_domain_Ud}.
We extend it on $\X^d$ by $0$. Without loss of generality, we can assume that 
\begin{equation*}
   \int_{\X^d} u_{d+1}^2 d\si_d=1
\end{equation*}
and by definition 
\begin{equation*}
   \int_{\X^d} |\nabla u_{d+1} |^2 d\si_d= \la_1(d+1).
\end{equation*}
Plugging this estimates in the inequality \eqref{sob_ineq} and using H\"older inequality, we get that 
\begin{equation*}
   \begin{aligned}
1  = \int_{\X^d} u_{d+1}^2 d\si_d
& \leq \vol_{\si_d}(U_{d+1})^{\frac{2}{d}}  \left( \int_{\X^d} |u_d|^{\frac{2d}{d-2}} d\si_d \right)^{\frac{d-2}{d}} \\
 & \leq \vol_{\si_d}(U_{d+1})^{\frac{2}{d}}  \left(\frac{ 4 \la_1(d+1)}{d(d-2) \om_{d}^{\frac{2}{d} }} + \frac{1}{  \om_{d}^{\frac{2}{d} }}\right).
\end{aligned}
\end{equation*}
We now notice that $ \vol_{\si_d}(U_{d+1})  = \left(1 -\frac{1}{2^d} \right)\om_d$. If we set $x_d=  \left(1 -\frac{1}{2^d} \right)^{\frac{2}{d}}$,  the last inequality then reads as 
\begin{equation*}
   \la_1(d+1) \geq \frac{d(d-2) (1 - x_d)}{  4 x_d}.
\end{equation*}
We finally write that when $d\to \infty$, $1- x_d \sim   \frac{2}{d} \frac{1}{2^d}$ and we obtain  
\begin{equation*}
   \la_1(d+1) \geq c  \frac{d}{2^d}
\end{equation*}
for some $c>0$, 
which gives the desired lower bound in Theorem~\ref{thm:main_estimate}.

\subsection{The upper bound}
Assume that $z_d:= \frac{2^d \la_1(d)}{d^3}$ is not bounded. Up to extracting a subsequence, we can assume that 
 \begin{equation} \label{hyp_cont}
 \lim_{d \to \infty} z_d=  \infty.
 \end{equation} 
 Choose a sequence $\al_d$ tending to $0$ such that 
 \begin{equation} 
 \label{cons_ald}
   z_d \gg \al_d^{-2}
 \end{equation}
and set $a_d= \al_d d^{-\frac{3}{2}}$. Let $\theta_d: \R \to \R$ be a smooth function satisfying $0 \leq \theta_d \leq 1$, $\theta_d(t) = 0$ if $t \geq 0$,  $\theta_d(t) = 1$ if $t \leq -a_d$ and $|\theta'| \leq 2 a_d^{-1}$.
For all $x \in \R^d$, set $\eta_d(x) =\max\{\theta_d(x_1), \ldots, \theta_d(x_d)\}$. Then the restriction of $\eta_d$ (still denoted by $\eta_d$) to $\X^{d-1}$ has the following properties: 
\begin{itemize}
 \item the function  $\eta_d$ is supported in $U_d$ since when $x \not\in U_d$ then $x_i>0$ for all $i$ and thus $\theta_d(x_i)=0$ for all $i$;
 \item $\nabla \eta_d$ exists almost everywhere and satisfies $|\nabla_{\X^{d-1}}\eta_d| \leq |\nabla_{\R^d} \eta_d | \leq 2a_d^{-1}$. 
\end{itemize}
Introduce the domains 
\begin{equation*}
    \Sigma_1 := \bigl\{x \in \X^{d-1} | \eta_d(x) = 1 \bigr\} \quad \hbox{and}\quad  \Sigma_{01}:= \bigl\{x \in \X^{d-1} | 0< \eta_d(x) < 1 \bigr\}.
\end{equation*}

\begin{step} \label{step1}
 It holds that  $|\Sigma_{01}|_{d-1} \leq C\left( \frac{1}{2} \right)^{d}  \om_{d-1}$. 
\end{step}

\begin{proof}
Using our notation \eqref{eq:def_dom}, we notice that 
\begin{equation*}
   \Si_{01} \subset \Si_{d-1}(a_d)=\{x \in \X^{d-1} | \forall i, \;  x_i \geq -a_d\}.
\end{equation*}
Indeed, as soon as $x_i \leq -a_d$ for some $i$, it holds that $\eta_d(x)=1$. 
 Since $a_d = \al_d  d^{-\frac{3}{2}} \ll d^{-\frac{3}{2}}$, we get from Lemma \ref{volsi}  that  
\begin{equation} \label{sigma}
    |\Si_{d-1}(a_d)|_{d-1} \leq C\left( \frac{1}{2} \right)^{d}  \om_{d-1},
\end{equation}
which proves Step~\ref{step1}. 
\end{proof}

\begin{step} \label{step2}
There exists $C>0$ such that  $|\Sigma_1|_{d-1} \geq  C \om_{d-1}$. 
\end{step}

\begin{proof}
We first observe that 
\begin{equation*}
   \Si'_d(a_d):=\{x \in \X^{d-1} | x_d \leq -a_d \} \subset \Si_1,
\end{equation*}
since $\theta_d(x_d)=1$ as soon as $x_d \leq -a_d$,  and thus 
$\eta_d(x) =1$. With the same notation as in Step~\ref{step1}, we thus have   
\begin{equation*}
     \X^{d-1} \setminus \Si_{d-1}(a_d) \subset \Si_1.
\end{equation*}
Since by \eqref{sigma}
\begin{equation*}
 | \Si_{d-1}(a_d)|_{d-1} \leq o(1) \om_{d-1},
\end{equation*}
we obtain that 
\begin{equation*}
  |\Si_1|_{d_1} \geq \om_{d-1} (1 - o(1)),
\end{equation*}
which implies Step \ref{step2}. 
\end{proof}

\begin{step} \label{conclusion}
 Conclusion. 
\end{step}

\begin{proof}
Recall that $\la_1(d)$ has the following variational characterization: 
\begin{equation*}
   \la_1(d) = \inf_{u \in H^1_0(U_d)\setminus \{0\}} Y(u),
\end{equation*}
where $H^1_0(\X^{d-1})$ is the set of $L^2(U_d)$-functions vanishing on $\partial U_d$ and whose derivatives of order $1$ are also $L^2$, and where 
\begin{equation*}
   Y(u) = \frac{ \int_{U_d} |\nabla u|^2 d\si_{d-1} }{\int_{U_d} u^2 d\si_{d-1}}.
\end{equation*}
In particular, it holds that $\la_1 (U_d) \leq Y(\eta_d)$.
Since the support of $\nabla \eta_d$ is included in $\Si_{01}$, it holds from Step \ref{step1} that  
\begin{equation*}
   \int_{U_d} |\nabla \eta_d|^2 d\si_{d-1}  \leq \|\nabla \eta_d\|^2_{L^{\infty} } |\Si_{01}|_{d-1}  \leq C a_d^{-2}  |\Si_{01}|_{d-1}  \leq C d^{3} \al_d^{-2}    \left(\frac{1}{2} \right)^{d}   \om_{d-1}, 
\end{equation*}
while by Step~\ref{step2}  
\begin{equation*}
   \int_{U_d} \eta_d^2 d\si_{d-1} \geq |\Si_1|_{d-1} \geq C \om_{d-1} .
\end{equation*}
This implies that $ \la_1(d) \leq Y(\eta_d) \leq C d^3  \al_d^{-2}   \left(\frac{1}{2} \right)^{d} $,  which contradicts \eqref{hyp_cont} and \eqref{cons_ald}. We thus obtain that 
\begin{equation*}
   \la_1(d) \leq C d^3  \left(\frac{1}{2} \right)^{d},
\end{equation*}
which proves the desired upper bound. 
\end{proof}

\subsection{Proof of Lemma~\ref{volsi}} \label{prooflem}

Introduce for any $d\geq 2$ and $0\leq k\leq d$ the set 
\begin{equation*}
    V_{k,d} (a) = \bigl([-a,1]^k \times [0,1]^{d-k}\bigr)\cap \mathbb S^{d-1}.
\end{equation*}
The objective of Lemma~\ref{volsi} is to estimate the volume of $V_{d,d}$. We shall decompose the proof into three steps. 

\medskip

In the first step, we prove that for all $a\in[0,1]$, all $d\geq 1$ and $0\leq k\leq d$,
\begin{equation}
    \label{eq:upper_bound_V_kd}
    \bigl\vert V_{k+1,d+1}(a)\bigr\vert_d \leq \bigl\vert V_{k,d+1}(a)\bigr\vert_d+a \bigl\vert V_{k,d}(\tfrac{a}{\sqrt{1-a^2}})\bigr\vert_{d-1}.
\end{equation}
Let us write $V_{k+1,d+1}(a)=V_+\cup V_-$, where
\begin{equation*}
    \left\{\begin{array}{lcl}
    V_+ &=& V_{k+1,d+1}(a) \cap \{x_{k+1}\geq0\},\smallskip\\
    V_- &=& V_{k+1,d+1}(a) \cap \{x_{k+1}<0\}.
    \end{array}
    \right.
\end{equation*}
First, let us note that $V_+=V_{k,d+1}(a)$. Moreover, write $\bigl\vert V_-\bigr\vert_{d} = \int_0^a \bigl\vert H_s\bigr\vert_{d-1} ds$, with 
\begin{equation*}
   H_s =V_{k+1,d+1}(a) \cap \{x_{k+1}=s\}.
\end{equation*}
We now set $y_1=\frac{x_1}{\sqrt{1-s^2}},\ldots,y_k=\frac{x_k}{\sqrt{1-s^2}}$ and $y_{k+1}=\frac{x_{k+2}}{\sqrt{1-s^2}},\ldots,y_d=\frac{x_{d+1}}{\sqrt{1-s^2}}$. With this notation, $(x_1,\ldots ,x_{d+1})\in H_s$ if and only if $x_{k+1}=s$ and $(y_1,\ldots, y_d)\in V_{k,d}\bigl(\frac{a}{\sqrt{1-s^2}}\bigr)$. We thus have
\begin{equation*}
    \bigl\vert H_s\bigr\vert_{d-1}\leq (1-s^2)^{\frac{d-1}{2}} \bigl\vert V_{k,d}\bigl(\tfrac{a}{\sqrt{1-s^2}}\bigr)\bigr\vert_{d-1}\leq (1-s^2)^{\frac{d-1}{2}} \bigl\vert V_{k,d}\bigl(\tfrac{a}{\sqrt{1-a^2}}\bigr)\bigr\vert_{d-1} \leq \bigl\vert V_{k,d}\bigl(\tfrac{a}{\sqrt{1-a^2}}\bigr)\bigr\vert_{d-1}.
\end{equation*}
We deduce that $\bigl\vert V_-\bigr\vert_{d}\leq a\bigl\vert V_{k,d}\bigl(\tfrac{a}{\sqrt{1-a^2}}\bigr)\bigr\vert_{d-1}$. The proof of \eqref{eq:upper_bound_V_kd} is complete.

\medskip

In the second step of the proof, we let $D\geq 2$ and prove the following property by induction on $d\in\{2,\ldots ,D\}$: there exist a quantity $\epsilon_D=O(\frac{1}{\sqrt{D}})$
such that
\begin{equation}
\label{eq:induction_H_d}
   \forall 0\leq k\leq d, \quad \forall a\leq \frac{\epsilon_D}{\sqrt{d-(d-k)\epsilon_D^2}},\quad
   \bigl\vert V_{k,d}(a)\bigr\vert_{d-1} \leq 4 \left(\frac{1}{2}+\epsilon_D\right)^{k} \frac{\om_{d-1}}{2^{d-k}}.
\end{equation}   
The property \eqref{eq:induction_H_d} is clearly valid for $d=2$.  Let us now assume that \eqref{eq:induction_H_d} holds at some $d$, and prove it at $d+1$. We will show the latter property by induction on $k$. It is clearly satisfied at $k=0$, using the following: the value of $V_{k,d+1}$ is easily calculated for $k=0$, at which
\begin{equation*}
   \bigl\vert V_{0,d+1}(a)\bigr\vert_{d} =\bigl\vert [0,1]^{d+1} \cap \mathbb S^{d}\bigr\vert_{d} = \frac{\omega_{d}}{2^{d+1}}.
\end{equation*}
%Recall also that $\sup_{d\geq 1}\omega_{d-1}=\omega_{6}<\infty$. 
Assuming it holds at a given $k$, we want to prove that
\begin{equation}
\label{eq:induction_F_k}
 \forall a\leq \frac{\epsilon_D}{\sqrt{d+1-(d-k)\epsilon_D^2}},\quad
   \bigl\vert V_{k+1,d+1}(a)\bigr\vert_d \leq 4 \left(\frac{1}{2}+\epsilon_D\right)^{k+1} \frac{\om_{d}}{2^{d-k}}.
\end{equation}  
To that aim, we use the upper bound \eqref{eq:upper_bound_V_kd}, which is a sum of two terms:
\begin{itemize}
    \item to get an upper bound of the term $\bigl\vert V_{k,d+1}(a)\bigr\vert_d$, we use the induction hypothesis at $k$;
    \item to obtain an upper bound of $\bigl\vert V_{k,d}(\tfrac{a}{\sqrt{1-a^2}})\bigr\vert_{d-1}$, we use the induction hypothesis in the index~$d$.
\end{itemize}
To apply the second upper bound, we need to verify that $\frac{a}{\sqrt{1-a^2}}$ satisfies the inequality in \eqref{eq:induction_H_d}. Starting from the inequality on $a$ in \eqref{eq:induction_F_k}, we easily obtain
\begin{equation*}
    \frac{a}{\sqrt{1-a^2}}\leq \frac{\epsilon_D}{\sqrt{d-(d-k)\epsilon_D^2+(1-\epsilon_D^2)}}\leq \frac{\epsilon_D}{\sqrt{d-(d-k)\epsilon_D^2}}. 
\end{equation*}
As a consequence, using \eqref{eq:upper_bound_V_kd}, \eqref{eq:induction_H_d} and \eqref{eq:induction_F_k}, we obtain
\begin{equation*}
    \bigl\vert V_{k+1,d+1}(a)\bigr\vert \leq 4 \left(\frac{1}{2}+\epsilon_D\right)^{k} \frac{1}{2^{d-k}}\omega_d\left(\frac{1}{2}+a\frac{\omega_{d-1}}{\omega_d}\right).
\end{equation*}
Using the fact that for all $d\geq 1$, $\frac{\omega_{d-1}}{\omega_d}\leq \sqrt{d}$, we easily obtain
\begin{equation*}
    a\frac{\omega_{d-1}}{\omega_d} \leq \frac{\sqrt{d}\epsilon_D}{\sqrt{d+1-d\epsilon_D^2}}\leq \epsilon_D,
\end{equation*}
assuming in the latter inequality that $d\epsilon_D^2\in(0,1]$.

\medskip

In the last part of the proof, we apply \eqref{eq:induction_H_d} at $d=D$ and deduce that for all 
$a\leq\frac{\epsilon_D}{\sqrt{D}}$,
\begin{equation*}
    \bigl\vert V_{D,D}(a)\bigr\vert_{D-1}\leq C\left(\frac{1}{2}+\epsilon_D\right)^D \om_{D-1}.
\end{equation*}
As soon as $\frac{\epsilon_D}{D}\to 0$, we have $\bigl(\frac{1}{2}+\epsilon_D\bigr)^D\leq C \frac{1}{2^D}$ for some $C>0$. The proof of Lemma~\ref{volsi} is completed.

\section*{Acknowledgments}
We would like to thank Michel Bauer and Françoise Cornu for suggesting the problem. KR would like to thank Nicolas Raymond and Grégory Schehr for interesting discussions. KR thanks the VIASM (Hanoï, Vietnam) for their hospitality and wonderful working conditions.

\end{document}